\documentclass[10pt]{amsart}
\usepackage{amsmath}
\usepackage{amssymb}
\usepackage{amsfonts}
\usepackage{amsthm}
\usepackage{enumerate}
\usepackage[mathscr]{eucal}

\theoremstyle{plain}
\newtheorem{theorem}{Theorem}[section]
\newtheorem{lemma}[theorem]{Lemma}
\newtheorem{corollary}[theorem]{Corollary}
\newtheorem{proposition}[theorem]{Proposition}

\numberwithin{equation}{section}

\theoremstyle{plain}

\numberwithin{equation}{section}

\theoremstyle{remark}
   \newtheorem{remark}[theorem]{Remark}

\def\inn#1#2{\langle#1,#2\rangle}

\def\eps {{\varepsilon}}

\def\vth {{\vartheta}}

\def\card{\text{card}}

\def\bbR{{\mathbb R}}

\def\lc{\lesssim}

\def\leps{{\,\lessapprox\, }}

\newcommand {\supp} {\mbox{supp }}

\newcommand {\SR} {{\Bbb R}}

\newcommand {\SZ} {{\Bbb Z}}

\renewcommand {\phi} {{\varphi}}

\newcommand {\al} {{\alpha}}
\newcommand {\dt} {{\delta}}
\newcommand {\Dt} {{\Delta}}
\newcommand {\e} {{\varepsilon}}
\newcommand {\ga} {{\gamma}}
\newcommand {\Ga} {{\Gamma}}
\newcommand {\la} {{\lambda}}
\newcommand {\La} {{\Lambda}}

\newcommand {\om} {{\omega}}

\renewcommand {\O} {{\Omega}}

\newcommand {\fS} {{\frak S}}

\newcommand {\cA} {{{\mathcal A}}}

\newcommand {\cD} {{\mathcal D}}

\newcommand {\cP} {{\mathcal P}}
\newcommand {\cQ} {{\mathcal Q}}

\newcommand {\cS} {{\mathcal S}}

\newcommand {\cT} {{\mathcal T}}

\newcommand {\cU} {{\mathcal U}}
\newcommand {\cW} {{\mathcal W}}

\newcommand {\tS} {{\widetilde S}}

\newcommand {\hf} {{\hat f}}

\newcommand {\hpsi} {{\widehat\psi}}

\renewcommand {\tfrac}[2]{{\textstyle\frac{#1}{#2}}}

\newcommand {\Scirc} {\raise.2ex\hbox{$\scriptstyle\circ$}}

\newcommand {\mand} {{\quad\mbox{and}\quad}}

\renewcommand {\mid} {{\,\,\,\colon\,\,\,}}

\newcommand{\Ba}[1]{\begin{array}{#1}}
\newcommand{\Ea}{\end{array}}
\newcommand{\Be}{\begin{equation}}
\newcommand{\Ee}{\end{equation}}
\newcommand{\Bea}{\begin{eqnarray}}
\newcommand{\Eea}{\end{eqnarray}}
\newcommand{\Beas}{\begin{eqnarray*}}
\newcommand{\Eeas}{\end{eqnarray*}}
\newcommand{\Benu}{\begin{enumerate}}
\newcommand{\Eenu}{\end{enumerate}}
\newcommand{\Bi}{\begin{itemize}}
\newcommand{\Ei}{\end{itemize}}

\newcommand{\BR}{\begin{Remark} \em}
\newcommand{\ER}{\end{Remark}}

\newcommand{\Proof}{\begin{proof}}
\newcommand{\ProofEnd}{\end{proof}    }

\newcommand{\BE}{\begin{example}}
\newcommand{\EE}{\end{example}}

\newcounter{remark}

\newcommand {\SRd} {{\SR^{d+1}}}

\begin{document}

\title
[Plate decompositions of cone multipliers] {On plate
decompositions of cone multipliers}


\author[G. Garrig\'os and A. Seeger]
{Gustavo Garrig\'os and Andreas Seeger}

\address{G. Garrig\'os \\
Dep. Matem\'aticas C-XV \\
Universidad Aut\'onoma de Madrid\\
28049 Madrid, Spain}
\email{gustavo.garrigos@uam.es}

\address{A. Seeger   \\
Department of Mathematics\\ University of Wisconsin-Madison\\Madison, WI 53706, USA}
\email{seeger@math.wisc.edu}

\begin{thanks}
{Some results of this  paper were  presented by G.G. at the Conference on Harmonic Analysis and
its Applications at Sapporo, 2005, and a preliminary version has
been included in the Hokkaido
University Technical Report series in Mathematics
(No. 103, December 2005).
G.G. thanks
Hokkaido University and Prof. Tachizawa for the invitation and support.
 Special thanks also to Aline Bonami for many suggestions and discussions
related to this topic.}
\end{thanks}\begin{thanks} {G.G. was  supported in part by the
European Commission, within the IHP Network ``HARP 2002-2006'',
contract number HPRN-CT-2001-00273-HARP, and also by
\emph{Programa Ram\'on y Cajal} and grant
``MTM2004-00678'', MCyT (Spain). A.S. was  supported in part by
NSF grant DMS 0200186.}
\end{thanks}

\begin{abstract}
An important inequality due to Wolff on
plate decompositions of cone multipliers
is known to have consequences for
sharp $L^p$ results on cone multipliers, local smoothing for the
wave equation, convolutions with radial kernels, Bergman
projections in tubes over cones, averages over finite type
curves in $\mathbb{R}^3$ and associated maximal functions.
We observe that  the range of $p$ in  Wolff's  inequality,
 for the conic and the spherical versions,
can be improved
by using bilinear restriction results.
We also use this inequality to
give some improved estimates on square functions
associated to decompositions
of  cone multipliers in low dimensions. This gives a new  $L^4$ bound
for the  cone multiplier operator in $\bbR^3$. \end{abstract}

 \maketitle

\section{Introduction}  \setcounter{equation}{0}

Let $\Ga=\{(\tau,\xi)\in\SR\times\SR^d\mid \tau=|\xi|\,\}$ denote
the forward light-cone in $\SRd$, $d\geq2$. For fixed $c>0$
and small $\dt>0$, we consider $\dt$-neighborhoods of the
truncated cone

    \[
\Ga_{\dt}(c)=\{(\tau,\xi)\in\SRd\mid 1\leq\tau\leq2\mand
\bigl|\tau-|\xi|\bigr|\leq c\dt\},
    \]
with the usual decomposition into plates subordinated to a
{\it $\sqrt\dt$-separated} sequence in the sphere $\{\om_k\}\subset
S^{d-1}$:

    \Be
\begin{gathered}
\Pi_k^{(\dt)}=\Bigl\{(\tau,\xi)\in\Ga_{\dt}(c)\mid
\bigl|\frac{\xi}{|\xi|}-\om_k\bigr|\leq c'\sqrt\dt\,\Bigr\};
\\
\text{dist} (\omega_k,\omega_{k'})\ge \sqrt \delta \quad \text { if } k\neq k'.
    \end{gathered}
\Ee
Let
\Be\label{alphap}
 \alpha(p):= d(\frac 12-\frac 1p)-\frac 12,
\Ee
the standard Bochner-Riesz critical index in $d$ dimensions.
Then  Wolff's inequality
is the assertion that
for all
$\e>0$
    \Be
\Big\|\sum_k f_k\Big\|_p\,\leq\,C_\e\,
\delta^{-\alpha(p)-\eps}
\Bigl(\sum_k\|f_k\|_p^p\Bigr)^{1/p}
    \label{Wp}\Ee
provided that
\Be \supp \widehat{f_k}\subset \Pi_k^{(\delta)}
\label{platesupp}
.\Ee
The power $\alpha(p)$
is optimal for each
$p$ (except perhaps for $\e>0$), and the inequality is conjectured
to hold for all $p>2+\frac{4}{d-1}$. In his fundamental work
\cite{W1}, Wolff
developed a method to prove such
inequalities for large values of $p$, and obtained a positive answer
for $d=2$ and $p>74$. Subsequently the method has been extended in the
paper by \L aba and Wolff \cite{LW} to higher dimensions. It is shown there that \eqref{Wp}  holds  for
$p>2+\frac{32}{3d-7}$ when $d\geq3$ and $p>2+\frac8{d-3}$ when
$d\geq4$.
In this paper we modify the weakest part of their proof
to obtain a better range of exponents in all dimensions (see
Table \ref{table1} below). The improvement relies on certain square
function
bounds  which follow from Wolff's bilinear Fourier extension theorem,
\cite{W2}.

\bigskip

\begin{table}[h]
  \centering
  \begin{tabular}{|c|l|l|l|} \hline
 Dimension & \cite{W1}, \cite{LW}     & Improvements               & Conjecture         \\
\hline
 $d=2$     & $p>74$            & $p>p_2:=63+1/3$                       & $p>6$              \\
\hline
 $d=3$     & $p>18$            & $p>p_3:=15$                           & $p>4$              \\
\hline
 $d=4$     & $p>8.4$           & $p>p_4:=7.28$                         & $p>\frac{10}3$     \\
\hline
 $d\geq5$  & $p>2+\frac8{d-3}$ & $p>p_d:=2+\frac8{d-3}(1-\frac1{d+1})$ & $p>2+ \frac4{d-1}$ \\
\hline
\end{tabular}
  \caption{Range of exponents for the validity of (\ref{Wp}) for light-cones.}\label{table1}
\end{table}


\begin{theorem}\label{th1} Let $d\geq2$ and $p_d$ as in Table
\ref{table1}. Then, under the assumption
\eqref{platesupp}
the  inequality (\ref{Wp}) holds for
 all $\e>0$
and all $p\ge p_d$.
\end{theorem}

\noindent {\it Remark:} Various further and   more technical
improvements on the range of Theorem \ref{th1} (and by implication on the results of Corollaries \ref{cor3} and \ref{corbergman} below)
 have been
obtained by the authors, and also by Wilhelm Schlag.
After the  initial submission of this article  these improvements
have been combined and  included in a joint paper (\cite{GSS}).

\medskip

A similar result
can be proved for decompositions of spheres  in $\mathbb R^d$.
We now  let
    \[
\cU_{\dt}(c)=\{\xi\in\mathbb R^d :
\big||\xi|-1\big|\leq c\dt\},
    \]
and consider the  decomposition into
plates subordinated to a $\sqrt\dt$-separated sequence in the sphere $\{\om_k\}\subset
S^{d-1}$,
    \[
B_k^{(\dt)}=\Bigl\{\xi\in\cU_{\dt}(c):
\bigl|\xi/|\xi|-\om_k\bigr|\leq c'\sqrt\dt\,\Bigr\}.
    \]

%

\begin{theorem} \label{thsphere}
The analogue of Wolff's inequality for the sphere,
    \Be
\Big\|\sum_k f_k\Big\|_p\,\leq\,C_\e\,
\delta^{-\alpha(p)-\eps}
\Bigl(\sum_k\|f_k\|_p^p\Bigr)^{1/p}, \quad \supp \widehat{f_k}\subset B_k^{(\delta)},
    \label{Wpsphere}\Ee
 holds for  $p\ge 2+\frac{8}{d-1}- \frac{4}{(d-1)d}$ and all
$\eps>0$.
\end{theorem}

Again \eqref{Wpsphere}
is conjectured to hold for the optimal range $p>2+4/(d-1)$.
It has been
 known to hold for  $p>2+8/(d-1)$;
this follows from a modification  of the argument in \cite{LW},
see also \cite{LP}. Note that in two dimensions the range is improved from previously $p>10$ to $p>8$.

\medskip

\noindent{\it Remark:} Theorem \ref{thsphere} may be extended to convex surfaces
with nonvanishing Gaussian curvature and similarly
Theorem \ref{th1} may be extended to cones with $d-1$ positive principal
curvatures. This can  be  achieved by using
scaling and induction on scales
 arguments such as in \S2 of \cite{PS}, see also the article by \L aba and Pramanik \cite{LP} for related results.

\medskip

We  proceed to list some of the known implications of Theorem \ref{th1}.



\begin{corollary}\label{cor3}
Let $d\geq2$ and $p_d$ be  as in  Table \ref{table1}. Then \Benu

\item [(i)]
For all $p>p_d$, $\al>\frac{d-1}2-\frac dp$, we have
\begin{equation} \label{lsm}\Big(\int_1^2 \big\|e^{it\sqrt{-\Dt}}f\big \|_{L^p(\SR^d)}^p dt\Big)^{1/p} \lesssim\|f\|_{L^p_\al(\SR^d)}.
\end{equation}

\item [(ii)] For all $p\in (p_d,\infty)$,  $\al>\frac{d-1}2-\frac dp$
the Fourier multiplier
    \begin{equation}
\label{conemultiplier}
    m_\al(\tau,\xi)= (1-|\xi|^2/\tau^2)^\al_+\,
    \end{equation}
 defines
a bounded operator in $L^p(\SR^{d+1})$.

\item [(iii)] Let $K\in \cS'(\bbR^d)$ be radial, let
$\varphi\in C^\infty_0(\bbR^d\setminus \{0\})$
 so that  $\varphi$ is radial and
not identically zero, and let
 $\e>0$. Let $K_t= \mathcal F^{-1} [\varphi \widehat K(t\cdot)]$.
Then for all Schwartz functions $f$ and $1<r<\tfrac{p_d}{p_d-1}$
$$\|K*f\|_{r}\le C_\eps \sup_{t>0}
\Big(\int|K_t(x)|^r(1+|x|)^\varepsilon dx
\Big)^{1/r} \|f\|_r. $$

\item [(iv)]  Let $\chi\in C^\infty_0 (\Bbb R)$ and $s\mapsto \gamma(s)\in \bbR^3$
a smooth curve satisfying
 $\sum_{j=1}^n
|\inn{\theta}{\gamma^{(j)}(s)}|\neq 0$ for every unit vector $\theta$
 and every $s\in \supp \chi$.
For $t>0$ define the convolution operator $A_t$
by
$$A_t f(x)=\int f(x-t\gamma(s)) \chi(s) ds.$$
Suppose that  $\max\{n, 32+2/3\}<p<\infty$ . Then $A_t$ maps
$L^p(\Bbb R^3)$ into
the $L^p$-Sobolev space
$L^{p}_{1/p}(\Bbb R^3)$.
 Moreover the  maximal function
$Mf=\sup_t|A_t f|$ defines a bounded operator on $L^p(\Bbb R^3)$.
\Eenu

\end{corollary}

Parts (i), (ii), (iii)  are standard  consequences of
Theorem  \ref{th1}; see  \cite{W1} for (i) and the local version of (ii).
The global version follows by results on dyadic decompositions of
multipliers and
$L^p$ Calder\'on-Zygmund theory (see \cite{ca} or \cite{se}).
The proof of Theorem 1.6 in \cite{MS} together
with these arguments can be used to deduce (iii) from Theorem \ref{th1}.
For (iv) see \cite{PS}.

Besides the connection to cone multipliers a major
motivation for this paper
is the relevance of
inequalities for plate decompositions
 for the boundedness properties of the Bergman projection  in tube domains
over full light cones, see
\cite{BBPR}, \cite{BBGR}. Denote by
$\Delta(Y)=y_0^2-|y'|^2$ the Lorentz form and consider the forward light cone
on which $\Delta$  is positive;
$$\La^{d+1}=\{Y=(y_0,y')\in \bbR\times\bbR^d: y_0^2-|y'|^2>0, y_0>0\}.$$
Let
$\cT^{d+1}\subset \mathbb{C}^{d+1}$ be the tube domain over  $\Lambda^{d+1}$, i.e.
$$\cT^{d+1}= \bbR^{d+1}+i \Lambda^{d+1}.$$
Let $w_\ga(Y)=\Delta(Y)^\ga$ and
 consider the weighted space $L^{p}(\cT^{d+1},w_\ga)$
with norm
%
$$\|F\|_{p,\ga}=
\Big(\iint_{\cT^{d+1}}|F(X+iY)|^p
\Delta^{\ga}(Y) \, dY dX  \Big)^{1/p}.
$$
Let $\cP_\ga$ be the orthogonal projection mapping the weighted space
$L^2(\cT^{d+1}, w_\ga)$
to its subspace $\cA^{p}_{\gamma}$ consisting of the  holomorphic functions.
Only the case $\ga>-1$ is interesting since $\cA^{p}_{\gamma}=\{0\}$ for $\gamma\le -1$.
We are interested in the $L^p$ boundedness properties of $\cP_\ga$.
For  $\ga>-1$ the operator $\cP_\ga$ can only be bounded
 on
$L^p(\cT^{d+1}, w_\ga)$ in the range
\Be\label{bergmanrange}
1+\frac{d-1}{2(\ga+d+1)}<p<1+\frac{2(\ga+d+1)}{d-1},
\Ee
see e.g. Theorem 4.3 in \cite{BBGNPR}, and \eqref{bergmanrange} is indeed the conjectured range for $L^p$ boundedness
(except for $d=2$ and
$\ga\in(-1,-1/2)$, in which case there are additional
counterexamples for $p\geq 8+4\ga$, see \cite{BBGR}).



\begin{corollary}\label{corbergman}
Let $d\geq2$ and $p_d$ as in  Table \ref{table1}. Then  for all
 $\ga\geq\frac{d-1}2(p_d-\frac{2(d+1)}{d-1})$,
the Bergman projection $\cP_\ga$ is a bounded operator in
$L^p(\cT^{d+1}, w_\gamma)$ in the sharp range \eqref{bergmanrange}.
\end{corollary}

Beyond Corollary \ref{corbergman} both Theorem \ref{th1} and Theorem
 \ref{cone2d} below  have implications for the range of
boundedness of the Bergman projector
$\cP_\ga$ in natural weighted mixed norm spaces.
For the derivation of Corollary \ref{corbergman} and further
discussion of mixed norm estimates we refer to  \cite{BBGR}
({\it cf.} in particular Proposition 5.5 and Corollaries 5.12 and 5.17).

\medskip

 Our approach to Theorem  \ref{th1} is
 based on bilinear methods, for which we
consider a closely related inequality: \Be\Big \|\sum_k
f_k\Big\|_p\,\leq\,C_\al\,\dt^{-\al}\,
\Bigl(\sum_k\|f_k\|_p^2\Bigr)^{1/2}.
    \label{W2}
\Ee One can conjecture the validity of (\ref{W2}) for all $\al>0$
and all $2<p<2+ \frac4{d-1}$, but for the moment no positive
result for any such $p$ seems to be known. The limiting point
$p=\frac{2(d+1)}{d-1}$ should be the hardest case, since by
interpolation and H\"older's inequality it
implies both (\ref{W2}) and (\ref{Wp}) in all the
conjectured ranges. This kind of inequality arises naturally in
the study of weighted mixed norm
inequalities for  the
 Bergman projection operator $\cP_\ga$, see \cite{BBGR}.

We shall deduce Theorem
\ref{th1} by using a stronger version of
\eqref{W2} for $p=2(d+3)/(d+1)$,  but with a power of $1/\dt$ which is
(probably) not optimal. Namely under the assumption
\eqref{platesupp}  we have
\Be
\Bigl \|\sum_k
f_k\Bigr\|_{\frac{2(d+3)}{d+1}}
\,\leq\,C_\e\,
\delta^{-\frac{d-1}{4(d+3)}-\e}
\,\Bigl\|\,\bigl(\sum_k|f_k|^2\bigr)^{1/2}
\Bigr\|_{\frac{2(d+3)}{d+1}},
\label{sq} \Ee
for all $\e>0$.
We prove this inequality  in \S\ref{two} using the bilinear approach
of Tao and Vargas \cite[$\S5$]{TV} and  the optimal bilinear
cone extension inequality of T. Wolff \cite{W2},  see
 Proposition \ref{sqlinear} below.  By
Minkowski's inequality
and interpolation (\ref{sq}) trivially implies non optimal
estimates for the inequality (\ref{W2}) for all $p\in(2,\infty)$
(see Corollary \ref{cor1} below). In \S \ref{three} we use these to
refine a part of Wolff's proof of (\ref{Wp}) and obtain
the new sharp estimates for large $p$ announced in Table \ref{table1}.
In \S\ref{four}  we improve on some of the  square function
results in low
dimensions; these yield in particular
the following  estimate
for the cone multiplier in $\Bbb R^{2+1}$.
\begin{theorem} \label{544thm}
Suppose $\alpha >
\frac{5}{44}\big(\frac{p_2-4}{p_2- \frac {41}{11}}\big)$.
Then the cone Fourier multiplier $m_\alpha$ defines a bounded operator on $L^4(\bbR^3)$ and the local smoothing result \eqref{lsm} holds in two dimensions.
\end{theorem}
This is a small improvement over the known range $\alpha> 5/44$ which
follows from a combination of  \cite{TV} and \cite{W2}.

%

\medskip

\noindent {\it Notation.} We shall use the notation $A\lc B$
if there is a  constant
(which may depend on $d$) so that $A\le CB$.
For families $(A_\delta,B_\delta)$, $\delta\le 1$ we use  $A_\delta\leps B_\delta$ if  for every $\eps\in (0,1) $ there is a
constant $C_\eps$ so that $A_\delta\le C_\eps \delta^{-\eps} B_\delta$ for $\delta<1$.

\section{The bilinear estimate} \label{two}

Following the approach by Tao and Vargas, we first establish an
equivalence between linear and bilinear versions of (\ref{sq}),
which is a higher-dimensional analogue of Lemma 5.2 in \cite{TV}.

\begin{lemma}
\label{bil1} Let $d\geq2$, and suppose that for some $p\in[2,\infty)$
and $\al>\max\{0, (d-1)(1/4- 1/p)\}$
\Be \label{sqb}
\Bigl\|\,\Bigl(\sum_{\om_k\in\O}
f_k\Bigr)\,\Bigl(\sum_{\om_{k'}\in\O'}
f_{k'}\Bigr)\,\Bigr\|_{p/2}\,\leq \,C\,\dt^{-2\al}\,
\Bigl\|\,\Bigl(\sum_{\om_k\in\O}|f_k|^2\Bigr)^{1/2}\,\Bigr\|_p
\,
\Bigl\|\,\Bigl(\sum_{\om_{k'}\in\O'}|f_{k'}|^2\Bigr)^{1/2}\,\Bigr\|_p,
\Ee holds for all $f_k\in \cS(\SR^{d+1})$ with
$\supp\widehat f_k\subset\Pi^{(\dt)}_k$,
all pairs of $1$-separated
subsets $\O,\O'\subset S^{d-1}$ and all $\dt\ll1$. Then we also have
 \Be \label{sql} \Big\|\sum_k f_k\Big\|_{p}\,\leq
\,C'\,\dt^{-\al}\,
\Bigl\|\,\Bigl(\sum_{k}|f_k|^2\Bigr)^{1/2}\,\Bigr\|_p,\quad
\supp\widehat {f_k}\subset\Pi^{(\dt)}_k. \Ee
\end{lemma}

We remark that
the  restriction on $\alpha$ for $p>4$ is  never  severe.
To see this we note that the condition $(d-1)(1/4-1/p)\le \alpha(p)/2$ holds iff  $d\ge 2$ and that
\eqref{sql} cannot hold with $\alpha<\alpha(p)/2$; this can be
proved using Knapp examples.


\Proof [Proof of Lemma \ref{bil1}]

Let $\Phi:Q\equiv[0,1]^{d-1}\to S^{d-1}$ be a smooth
parametrization of (a compact subset of) the sphere and let $\cD$
denote the set of all dyadic intervals $I\subset Q$ with $|I|\geq
\dt^{\frac{d-1}2}$. As in \cite[p. 971]{TVV}, we may consider a
Whitney decomposition
$Q\times Q
    \,=\,\biguplus_{I\sim J}\,I\times    J$,
where $I\sim J$ means:
\medskip

\noindent (i) $I,J\in\cD$ and $|I|=|J|$;
\medskip

\noindent (ii) If $|I|>\dt^{\frac{d-1}2}$, then $I$ and $J$ are
not adjacent but their parents are.
\medskip

\noindent (iii) If $|I|=\dt^{\frac{d-1}2}$, then $I,J$ have
adjacent or equal parents.
\medskip

\noindent For simplicity, we
assume (by splitting the sphere
into finitely many pieces) that all $\om_k\in \Phi(Q)$ and let
 $y_k=\Phi^{-1}(\omega_k) \in Q$. We also
denote $\cD_j=\{I\in\cD\mid |I|=2^{-j(d-1)}\}$. Then
    \[
\Bigl(\sum_{k} f_k\Bigr)^2\,= \,\sum_{y_k,y_{k'}\in Q}
f_k\,f_{k'}\, = \,
\sum_{\sqrt\dt\leq2^{-j}\leq1}\,\sum_{{I,J\in\cD_j}\atop{I\sim
J}}\, \bigl(\sum_{y_k\in I} f_k\bigr)\,\bigl(\sum_{y_{k'}\in
J} f_{k'}\bigr). \] To establish (\ref{sql}) we take
$L^{p/2}$-norms in the above expression and use Minkowski's inequality
in $j$, so that we reduce the problem to show, for each $j$ \Be
\Bigl\|\,\sum_{{I,J\in\cD_j}\atop{I\sim J}}\,
\Bigl(\sum_{y_k\in I} f_k\Bigr)\,\Bigl(\sum_{y_{k'}\in J}
f_{k'}\Bigr)\,\Bigr\|_{p/2}\,\lesssim\,(2^{2j}\delta)^{-2\alpha}
\max\{1, 2^{j(d-1)(1-4/p)}\}
\,
\Bigl\|\,\Bigl(\sum_{k}|f_k|^2\Bigr)^{1/2}\,\Bigr\|^2_p.\label{sqb2}
\Ee

Inequality \eqref{sqb2}  is trivial when $2^{-j}\approx \sqrt\dt$ since by assumption the
number of $y_k$'s in each $I$ is approximately constant. We
consider the general case $\sqrt\dt<2^{-j}\leq1$. By construction
we must have \Be\sum_{I\in\cD_j}\sum_{J\sim I}\chi_{I+J}\lesssim
1\label{og}.\Ee Indeed, if $c_I$ denotes the center of $I$, then
    \[
I+J\,\subset
\,\bigl(c_I+B_{c2^{-j}}\bigr)\,+\,\bigl(c_J+B_{c2^{-j}}\bigr)\,
\subset\,2c_I\,+\,B_{c'2^{-j}}.
    \]
Since for each $I$  there are at most $O(1)$ cubes $J$ with  $J\sim I$,
and since  the centers $c_I$ are $2^{-j}$ separated, (\ref{og}) follows
easily.

From (\ref{og}) it follows that the functions
$F_{I,J}=\bigl(\sum_{y_k\in I}
f_k\bigr)\,\bigl(\sum_{y_{k'}\in J} f_{k'}\bigr)$ have pairwise
(almost) disjoint spectra when $I\sim J\in\cD_j$.
We may conclude  by orthogonality and standard interpolation arguments
\Be \Big\|\sum_{I\sim
J\in\cD_j} \,F_{I,J}\,\Big\|_{p/2}\,\lesssim\,
\max\{1, 2^{j(d-1)(1-4/p)}\}
\Bigl(\,\sum_{I\sim
J\in\cD_j} \,\|F_{I,J}\|^{p/2}_{p/2}\,\Bigr)^{2/p}.
\label{og1} \Ee (the case $p/2=2$ follows by orthogonality and the cases
$p/2=1$ and $p/2=\infty$ are trivial; see e.g. Lemma 7.1 in \cite{TV}). Next, we wish
to use the bilinear assumption (\ref{sqb}) to estimate
$\|F_{I,J}\|_{p/2}$. This can only be used directly when $2^j\approx 1$,
since dist$ (I,J)\sim1$. For other $j$'s we must use Lorentz
transformations to rescale the problem. To do this, let
$\{\eta_1,\ldots,\eta_d\}$ be an orthonormal basis of $\SR^d$
with $\eta_1$ being the center of $\Phi(I)$. Then we define $L\in
SO(1,d)$ acting on a basis of $\SR^{d+1}$ by
    \[
L(1,\eta_1)=(1,\eta_1), \quad
L(-1,\eta_1)=\tfrac\sigma\dt\,(-1,\eta_1)\mand
L(0,\eta_\ell)=\sqrt{\tfrac\sigma\dt}\,(0,\eta_\ell),\;\ell=2,...,d,
    \]
where we choose $\sigma=2^{2j}\dt$ (so that $\delta<\sigma<1$).
The functions
$f_k\circ L$ have now spectrum in
(perhaps a multiple) of the plates $\Pi^{(\sigma)}_k$
corresponding to the $\sqrt\sigma$-separated centers
$\{L(1,\om_k)\}$. Moreover, by the choice of $\sigma$,
the plates
corresponding to $y_k\in I$ and $y_{k'}\in J$ are $c$-separated,
and therefore after a change of variables we can apply
(\ref{sqb}) at scale $\sigma$ to obtain
    \begin{align*}
\|F_{I,J}\|_{p/2}\,&=\,\Bigl\|\,\Bigl(\sum_{y_k\in I}
f_k\Bigr)\,\Bigl(\sum_{y_{k'}\in J}
f_{k'}\Bigr)\,\Bigr\|_{p/2}\,\\&\lesssim
\,(2^{2j}\dt)^{-2\alpha}\,
 \Bigl\|\,\Bigl(\sum_{y_k\in I
}|f_k|^2\Bigr)^{1/2}\,\Bigr\|_p \,
\Bigl\|\,\Bigl(\sum_{y_{k'}\in J
}|f_{k'}|^2\Bigr)^{1/2}\,\Bigr\|_p,
\end{align*}
and then also
\begin{align*}
\Bigl(\,\sum_{I\sim
J\in\cD_j} \,&\|F_{I,J}\|^{p/2}_{p/2}\,\Bigr)^{2/p}
\lc
(2^{2j}\dt)^{-2\alpha}\,
\Bigl[\,\sum_{I\sim
J\in\cD_j} \,\Bigl\|\Bigl(\sum_{y_k\in I
}|f_k|^2\Bigr)^{1/2}\Bigr\|^\frac p2_p \,\,
\Bigl\|\Bigl(\sum_{y_{k'}\in J
}|f_{k'}|^2\Bigr)^{1/2}\Bigr\|^\frac p2_p\,\Bigr]^\frac2p
\\ &\lc
(2^{2j}\dt)^{-2\alpha}\,
\Bigl[\,\int
\Bigl(\sum_{I}
\sum_{y_k\in I }|f_k|^2\Bigr)^{p/2}
\Bigr]^{2/p}
\leq
(2^{2j}\dt)^{-2\alpha}\,
\Bigl\|\,\Bigl(\sum_{k}|f_k|^2\Bigr)^{1/2}\Bigr\|^2_p,
 \end{align*}
where in the second inequality  we have
used  $2ab\leq a^2+b^2$ followed by the imbedding $\ell^1\hookrightarrow\ell^\frac p2$.
Combining this with \eqref{og1}
we obtain
\Be \Big\|\sum_{I\sim
J\in\cD_j} \,F_{I,J}\,\Big\|_{p/2}\,\lesssim\,
(2^{2j}\delta)^{-2\alpha} \max\{1, 2^{j(d-1)(1-4/p)}
\}
\Bigl\|\Bigl(\sum_{k}|f_k|^2\Bigr)^{1/2}\Bigr\|^2_p,
\Ee
This  proves (\ref{sqb2}).
By our assumption on $\alpha$ we may sum in $j$
 and  the lemma follows.
\ProofEnd

We turn to the proof of (a generalization of) the square function estimate
\eqref{sq}.
We shall use the following
statement of Wolff's Fourier extension  theorem.


\medskip
\noindent {\bf Wolff's bilinear estimate.}  \cite[p.
680]{W2}. \emph{Let $p\ge \frac{d+3}{d+1}$, $\e>0$ and let $E,E'$ be
1-separated subsets of $\Ga_{1/N}$. Then, for all smooth $f$ and
$g$ supported in $E$ and $E'$, and all $N$-cubes $Q$, we have}
\begin{equation}
\label{Wbil}
\bigl\|\widehat f\,\widehat g\bigr\|_{L^p(Q)}
\,\leq\,C_\e\,N^{-1+\e}\,\|f\|_2\,\|g\|_2.
   \end{equation}
\medskip

Denote by
$\cQ\equiv \cQ(\delta^{-1/2})$ a tiling of
$\bbR^{d+1}$ with cubes $Q$ of disjoint interior and sidelength
$\delta^{-1/2}$,
with
centers $c_Q$ in $\dt^{-\frac12}\SZ^{d+1}$.

\begin{proposition}\label{propbilinear}
Let $d\geq2$, and suppose that
$\supp \widehat{f_k}\subset \Pi_k^{(\delta)}$,
$\supp \widehat{g_k}\subset \Pi_k^{(\delta)}$
and let   $\O,\O'\subset S^{d-1}$ be $1$-separated subsets.
Suppose  $\frac{2(d+3)}{d+1}\le q\le p\le \infty$
and let
\begin{equation} \label{mu} \mu(p)=
\frac {d}{4}-\frac {d+1}{2p}.
\end{equation}
Then, for all $\e>0$
\begin{multline}
\Big(\sum_{Q\in \cQ(\delta^{-1/2})}
    \Bigl\|\,\bigl(\sum_{\om_k\in\O}
f_k\bigr)\,\bigl(\sum_{\om_{k'}\in\O'}
g_{k'}\bigr)\,\Bigr\|_{L^{q/2}(Q)}^{p/2}\Big)^{2/p}\\
\,\lesssim\,\,
\delta^{-2\mu(p)-\eps}
\Bigl\|\,\Bigl(\sum_{\om_k\in\O}|f_k|^2\Bigr)^{1/2}\,\Bigr\|_p
\,
\Bigl\|\,\Bigl(\sum_{\om_{k'}\in\O'}|g_{k'}|^2\Bigr)^{1/2}\,\Bigr\|_p.
\label{ss2}\
\end{multline}
\end{proposition}

\begin{proof}


Let $\psi\in\cS(\SR^{d+1})$ be so that $\supp\hpsi\subset
B_{1/10}$ and $\psi(x)>1$ if $|x_i|\le 2$, $i=1,\dots, d+1$; then
 $\sum_{n\in\SZ^{d+1}}\psi(\cdot+n)^2\approx 1$.
Let
$\psi_Q=\psi(\sqrt\dt(\cdot-c_Q))$, so that
$\sum_{Q}\psi_Q^2\approx 1$. We write
    \[
F^Q=\bigl(\sum_{\om_k\in\O} f_k\bigr)\,\psi_Q\mand
G^Q=\bigl(\sum_{\om_{k'}\in\O'} g_{k'}\bigr)\,\psi_Q,\] so that
the supports of
$\widehat{F^Q}$ and $\widehat{G^Q}$ are 1-separated sets in
$\Ga_{\sqrt\dt}$. Thus, 
%
we can use Wolff's estimate  \eqref{Wbil} with $N=\dt^{-1/2}$ 
to obtain \Be \Big\|\Big(\sum_{\om_k\in\O}
f_k\Big)\,\Big(\sum_{\om_{k'}\in\O'} g_{k'}\Big)\Big\|_{L^{q/2}(Q)}\,
\lesssim\,\big\| F^Q\,G^Q\big\|_{L^{q/2}(Q)} 
\,\leps \dt^{1/2}\,
\big\|\widehat{F^Q}\big\|_2\,\big\|\widehat{G^Q}\big\|_2.\label{ss1} \Ee Now, by almost
orthogonality we can write

    \[\big\|\widehat{F^Q}\big\|^2_2\approx
\sum_k\big\|\widehat{f_k}*\widehat{\psi_Q}\big\|_2^2=
    \Big\|\Big(\sum_k|f_k|^2\Big)^{1/2}\,\psi_Q\Big\|^2_2,
    \]
and similarly for $G^Q$. We write
$S_\O=(\sum_{\om_k\in\O}|f_k|^2)^{1/2}$,
$\tS_{\O'}=(\sum_{\om_k\in\O'}|g_k|^2)^{1/2}$,
raise (\ref{ss1}) to
the power $p/2$ and sum in $Q$.
Thus
\begin{equation*}
\Big(\sum_Q \Bigl\|\,\bigl(\sum_{\om_k\in\O}
f_k\bigr)\,\bigl(\sum_{\om_{k'}\in\O'} g_{k'}\bigr)\,
\Bigr\|_{L^{q/2}(Q)}^{p/2}\Big)^{2/p}
\leps  \sqrt\delta \,
\Bigl(\sum_Q\,\bigl\|
S_\O\psi_Q\bigr\|_2^{p/2}\big\|\tS_{\O'}\psi_Q\big\|_2^{p/2}\Big)^{2/p}
\end{equation*} and by the Cauchy-Schwarz and H\"older inequalities the right hand side is

 \begin{align*}
&\lc \sqrt\delta
\Bigl(\sum_Q\,\bigl\|
S_\O\psi_Q\bigr\|_2^{p}
\Bigr)^{1/p}
\Bigl(\sum_Q\,\bigl\|\tS_{\O'}\psi_Q\big\|_2^{p}\Big)^{1/p}
\\
\;& \lc
\sqrt\delta
\Bigl(\sum_Q\,\bigl\|
S_\O\psi_Q\bigr\|_p^{p}
|Q|^{-1+p/2}
\Bigr)^{1/p}
\Bigl(\sum_Q\,\bigl\|\tS_{\O'}\psi_Q\big\|_p^{p} |Q|^{-1+p/2}\Big)^{1/p}
\\
& \lc  \dt^{\frac 12-(d+1)(\frac 12 -\frac 1p)}\,
\bigl\|S_\O\bigr\|_p\,\bigl\|\tS_{\O'}\bigr\|_p
\end{align*}
which yields the assertion.
\end{proof} 

We combine    Proposition  \ref{propbilinear} for $q=p$ and
Lemma \ref{bil1}
 to obtain
\begin{proposition}\label{sqlinear}
Let $d\geq2$, let $\mu(p)$ be as in \eqref{mu} and suppose that
$p\ge \frac{2(d+3)}{d+1}$.
Then, for all $\e>0$
 \Be
\Bigl \|\sum_k
f_k\Bigr\|_{p}\,\leq\,C_\e\,
\delta^{-\mu(p)-\e}
\,\Bigl\|\,\bigl(\sum_k|f_k|^2\bigr)^{1/2}\,\Bigr\|_{p}
\quad \text{ if }\supp \widehat{f_k}\subset \Pi_k^{(\delta)}.
\label{squarefctstandard}
\Ee
\end{proposition}

We may apply
Minkowski's inequality on the right hand side of \eqref{squarefctstandard}
and  obtain
\eqref{W2} for the limiting case $p=2(d+3)/(d+1)$
and $\alpha>\mu(p)= d/4-(d+1)/2p$.
 It
turns out this is all that is needed to obtain the claimed
improvements in Theorem \ref{th1}. The resulting  inequality can also be
interpolated with the trivial estimates for $L^2$ and $L^\infty$
to give:

\begin{corollary}
\label{cor1} The inequality (\ref{W2}) holds for all
$\al>\frac{d-1}4(\frac12-\frac1p)$, when  $2\le p\le \frac{2(d+3)}{d+1}$
and for all
$\al>\frac{d-1}4(1-\frac{2(d+2)}{p(d+1)})$ when
$\frac{2(d+3)}{d+1}
\le p\le \infty$.
\end{corollary}


\section{An improvement of Wolff's estimate}\label{three}
\setcounter{equation}{0}

We turn to Theorem \ref{th1}. The proof in \cite{W1, LW} for
inequality (\ref{Wp}) is based on a subtle localization procedure,
induction on scales and certain combinatorial arguments. Here we
only discuss  the  modifications leading to the claimed
improvements based on
Proposition \ref{sqlinear}. A more self-contained exposition with further improvements can  be  found in \cite{GSS}.

 For simplicity, when $\dt$ is fixed (and small)
we use  the notation $A\lessapprox B$ to indicate the inequality
 $A\leq C_\e\,\delta^{-\e} B$ for
all $\e>0$. Recall that the number of plates $\Pi_k^{(\dt)}$
covering $\Ga_\dt$ is approximately $\delta^{-\frac{d-1}2}$. Also,
throughout this section we fix $q(d)=2(d+3)/(d+1)$.

Due to various reductions (see \cite[$\S3$]{LW}), it is enough to
show that, for all $f_k$ with $\supp\hf_k\subset\Pi_k^{(\delta)}$ and
$\|f_k\|_\infty\leq1$, and for all $\la>0$ we have \Be
\bigl|\{|{\textstyle\sum_k f_k}|>\la\}\bigr|\,\lessapprox
\la^{-p}\,\delta^{d-\frac{(d-1)p}2}\,\|f\|_2^2
\label{weak}
\Ee
where $f=\sum_k f_k$.
In \cite{W1,LW} it is observed that, by Chebyshev's inequality,
this property  trivially holds for small enough $\la$; namely for
all $\la\leq \dt^{-\frac{d-1}2+\frac1{p-2}}$. We
use \eqref{W2} to enlarge this range of $\la$.

\begin{lemma}
Let $q=q(d)=2(d+3)/(d+1)$. Then, inequality (\ref{weak}) holds for all
\begin{equation}\label{larange}
\la\leq \delta^{-\frac{d-1}2+\frac {q}{4(p-q)}}.
\end{equation}
\end{lemma}

\Proof Let 
$\beta=\frac{d-1}{4(d+3)}$. By
Chebyshev's inequality and \eqref{W2}, we have
$$
\bigl|\{|f|>\la\}\bigr|  \leq
\la^{-q}\,\|f\|^q_q\;\lessapprox\;\delta^{-q\beta}\,\la^{-q}\,\bigl(\textstyle\sum_k\|f_k\|_q^2\bigr)^{q/2}
$$
and estimate
$$
\bigl(\textstyle\sum_k\|f_k\|_q^2\bigr)^{q/2}
\lc \dt^{-\frac{d-1}2\frac{q/2}{(q/2)'}}
 \textstyle\sum_k\|f_k\|_q^q \,
\\
\lc \dt^{-\frac{d-1}2(\frac
q2-1)}\,\textstyle\sum_k\|f_k\|_2^2\,
\sup_k\|f_k\|_\infty^{q-2}.$$
Since by assumption
$\|f_k\|_\infty\leq1$ and by almost orthogonality
$\sum_k\|f_k\|_2^2\approx\|f\|_2^2$, it suffices to show that in
the desired range of $\la$ we have $\dt^{-q\beta-\frac{d-1}2(\frac q2-1)}
\la^{-q}\leq\dt^{-\frac{(d-1)p}2-d}\la^{-p}$ which is equivalent to
\eqref{larange}.
\ProofEnd

At this point one can proceed exactly as in the proof of
Proposition 3.2 of \cite{LW} (or p. 1277 in \cite{W1}, when
$d=2$). The desired gain comes from using $\la\geq
\dt^{-\frac{d-1}2+\frac {q(d)}{4(p-q(d))}}$ (rather than $\la\leq
\dt^{-\frac{d-1}2+\frac1{p-2}}$) in step (54) of \cite{LW} (or (68)
of \cite{W1}).

For completeness, we shall briefly sketch this procedure here,
referring always to the notation in \cite{LW}. Localizing with
$\sqrt{N}$-cubes $\Dt$ as in Lemma 6.1 of \cite{LW}, one can find
a collection of functions $\{f_\Dt\}$ with spectrum in
$\Ga_{\sqrt\dt}$ and a number
\Be\la_*\in (\la \delta^{\frac{d-1}4+\e},
c\delta^{-\frac{d-1}4})
\label{lastar}
\Ee so that
$$ \bigl|\{|f|>\la\}\bigr| \,
\lessapprox \,\sum_\Dt \bigl|\{|f_\Dt|>\la_*\}\bigr|
$$ and
\Be
\card\big(\cP(f_\Dt)\big)\lessapprox
\la_*^2\la^{-2}\,\dt^{-\frac{3d-1}{4}}.\label{aux1}\Ee Here
$\cP(f_\Dt)$ refers to the set of plates in the wave-packet
decomposition of $f_\Dt$. When the cardinality of this set is
``small'', a further localization argument and induction on scales
allows to conclude the theorem (see Lemmas 6.2 and 6.3 in
\cite{LW}).

In \cite{LW,W1}, the size of $\card(\cP(f_\Dt))$ which
ensures the validity of these arguments is controlled in three
different ways, each depending on a different combinatorial
estimate
\Be
\card\big(\cP(f_\Dt)\big)
\le c_\eps
\delta^\eps \la^2_*
\label{aux5},
\Ee
or
\Be\card \bigl(\cP(f_\Dt)\bigr)\le c_\eps
\dt^{\frac{3d-3}8 +\eps}\,\la^4_*,
\label{aux34}\Ee
or, in three  dimensions (i.e. $d=2$) only,
\Be\label{aux2d}
\card\bigl(\cP(f_\Dt)\bigr)\le  c_\eps\dt^{\frac{11}8+\eps}\,\la^9_* .
\Ee
the last  estimate  being by far the most difficult
(see Lemmas 5.2 and 5.3 in \cite{LW} and Lemma
3.2 in \cite{W1}).

Given the lower bound for $\la_*$ in \eqref{lastar}
and
\Be\label{opplarange}
\la\geq \delta^{-\frac{d-1}2+\frac {q}{4(p-q)}}
\Ee
and given   \eqref{aux1}
it remains to verify
the estimates \eqref{aux5} in the
claimed range $p>p_d$, $d\ge 5$, \eqref{aux34} for $p>p_d$, $d=3,4$
 and \eqref{aux2d} for $p>p_2$.

This is straightforward.
By \eqref{aux1}  and \eqref{opplarange} we have
$$\card \big(\cP(f_\Delta)\big) \lessapprox \delta^{-\eps} \la_*^2
\delta^{d-1-\frac{q(d)}{2(p-q(d))}} \delta^{-\frac{3d-1}4}
$$
which gives in the case $d\ge 4$ the assertion \eqref{aux5} if
$d-1-\frac{q(d)}{2(p-q(d))}-\frac{3d-1}4>0$ or, after a short computation
 $p>q(1+\frac2{d-3})=2+\frac8{d-3}\frac
d{d+1}$. This is the asserted range if  $d\ge 5$.

Next we examine the validity of the inequality \eqref{aux34} under condition \eqref{opplarange}.
We now have
\[ \card\bigl(\cP(f_\Dt)\bigr)\;
\;
\le C_\eps
\frac{\la_*^4
\delta^{-\frac{3d-1}4-\eps}}{\la^2_*\la^2}\;\leq\;\frac{\la_*^4
\dt^{-\frac{3d-1}4-\eps}}{\la^4\,\dt^{\frac{d-1}2+2\e}} \;\leq\;\frac{
\dt^{-\frac{5d-3}4-3\e}}{\dt^{-2(d-1)+\frac {q(d)}{p-q(d)}}}\,\la^4_*.\]

This quantity is $\lesssim \delta^\eps
\dt^{\frac{3d-3}8}\,\la^4_*$ if and only if
$
\tfrac{5d-3}4-2(d-1)+\tfrac {q(d)}{p-q(d)}+4\eps<-\tfrac{3d-3}8,
$
which yields the range
gives  $p>q(d)(1+\frac{8}{3d-7}).$
Notice that this inequality amounts to $p>7.28$ if $d=4$ and $p>15$ if
$d=3$ which is the assertion in those cases.

Finally we consider the case $d=2$ when $q(2)=10/3$.
By \eqref{aux1} we need to have
$\lambda_*^2\la^{-2}\dt^{-5/4-\eps}\le c_\eps \delta^{11/8}\la_*^9$ , i.e.
$\la^{-2}\delta^{-21/8-\eps}\le c_\eps\lambda_*^7$
provided that
$\la_*>\la\delta^{1/4+\eps}$. Thus taking the smallest possible $\lambda^*$
yields $\delta^{-35/8-10\eps}\le \la^9$ and this has to
hold for all $\lambda$ satisfying \eqref{opplarange}, i.e.
$\la\ge \delta^{-\frac 12 +\frac{q(2)}{4(p-q(2))}}$.
Taking the minimal $\la$ this is achieved if
$35/8-10\eps<9/2 -9q/(4p-4q)$ with $q=q(2)=10/3$. Solving in $p$ and letting $\eps\to 0$ yields the range
$p>19q(2)=63+1/3$.\qed


\begin{proof}[On the proof of Theorem \ref{thsphere}]
The proof  is similar to the proof of Theorem  \ref{th1}. Instead
of \eqref{sq} we use a square function inequality for the sphere
\Be \Bigl \|\sum_k
f_k\Bigr\|_q\,\leq\,C_\e\,\delta^{-\alpha(q)/2-\eps}
\,\Bigl\|\,\bigl(\sum_k|f_k|^2\bigr)^\frac12\,\Bigr\|_q, \quad
\supp \widehat{f_k}\subset B_k^{(\delta)}, \label{sqsphere} \Ee
with $\alpha(q)= d(1/2-1/q)-1/2$, and $q= 2(d+2)/d$. In two
dimensions this is an old observation by C. Fefferman (\cite{F}),
and holds for $q=4$ with $\eps=0$. In higher dimensions the
inequality \eqref{sqsphere} was proved by Bourgain \cite{Bo1} for
the range of the Stein-Tomas restriction theorem (i.e.
$q\ge2(d+1)/(d-1)$). For the larger range $q>2(d+2)/d$
 the proof of \eqref{sqsphere}  is rather analogous to the proof of
Proposition \ref{sqlinear}; one now uses   Tao's bilinear Fourier extension  inequality \cite{tao2}
(see also \cite{Lee} for
related  results). Unlike \eqref{squarefctstandard} in the
conic case  the inequality \eqref{sqsphere} is essentially optimal
for the given range $q\ge 2(d+2)/d$.
 We omit further details.
\end{proof}

\section{
More on square functions} \label{four}
We shall now discuss some improvements of the square function estimate in
Proposition \ref{sqlinear} in low dimensions; thus we seek  estimates
of the form

\Be
\Bigl \|\sum_k
f_k\Bigr\|_{p}\,\leq\,C_\e\,
\delta^{-\beta-\e}
\,\Bigl\|\Bigl(\sum_k|f_k|^2\Bigr)^{1/2}\,\Bigr\|_{p}, \quad
\supp \widehat{f_k}\subset \Pi_k^{(\delta)}.
\label{sqestimate}
 \Ee
for some $\beta<\mu(p)=d/4-(d+1)/2p$.
 We shall assume
 throughout this section  the following Wolff hypothesis, labeled
 hypothesis $\cW(w;d)$, and  aim to prove estimates of the form
 \eqref{sqestimate} conditional on
this hypothesis.

\medskip

\noindent {\bf Hypothesis} \emph{ $\mathbf{\boldsymbol{\mathcal W}\boldsymbol(\boldsymbol{w};\boldsymbol{d}\boldsymbol)}$.
For all $\delta\in (0,1)$ and all families
$\{h_k\}$ of functions satisfying $\supp{\widehat h_k}\subset
\Pi_k^{(\delta)}$,}
\begin{equation}
\label{hypothesisW}
\Big\|\sum_k h_k \Big\|_{w}\le C_\eps \delta^{-\alpha(w)-\eps}
\Big(\sum_k\|h_k\|_w^w\Big)^{1/w},
\end{equation}
\emph{where $\alpha(w)= d(1/2-1/w)-1/2$.}
{\it Cf.}
Table \ref{table1}.

\medskip

We note that in view of the embedding
$L^p(\ell^2)\subset L^p(\ell^p)$ the inequality
\eqref{hypothesisW} trivially implies \eqref{sqestimate} with $\beta=\alpha(p)$, for $w\le p<\infty$.
Another trivial observation is that
\eqref{sqestimate} holds  with $\beta \ge (d-1)/4$
in view of the  Cauchy-Schwarz inequality, as
 $\sum_k  |f_k(x)| \lc \delta^{-(d-1)/4}
(\sum_k|f_k(x)|^2)^{1/2}$ for every $x$.

The method for our improvement over the exponent
$\min\{ \mu(p), \frac{d-1}{4}\}$
will  be limited
to the case where
\Be\label{amuassumpt}
\alpha(p)<\min\{\mu(p),\tfrac{d-1}{4}\}
\Ee
which holds  if and only if
$p<\min\{\frac{2(d-1)}{d-2}, \frac{4d}{d-1}\}$.
We have the additional restriction 
$p> \frac{2(d+3)}{d+1}$ in Proposition \ref{sqlinear}. Summarizing
 we shall get an improvement which is limited to
$d=2,3,4$ and to the ranges
\begin{equation}
\label{threecases}
\begin{cases}
d=2, \quad& 10/3<p<\min\{8,w\},
\\
d=3, \quad& 3<p<4,
\\
d=4, \quad& 14/5<p<3.
\end{cases}
\end{equation}

We emphasize that square function estimates such as
\eqref{sqestimate} cannot {\it a priori} be interpolated when subject to the
Fourier support condition
\eqref{platesupp}.
We shall however start with a
preliminary result which is proved using an interpolation.

We let
$\phi_k$ be a bump function adapted to the plate $\Pi_k^{(\delta)}$
satisfying the natural estimates, so that
$\phi_k$ equals  $1$ on the plate, and is supported on the ``double plate''.
 Define the operator  $P_k$ by
\Be \label{pkdef} \widehat {P_k f}=\phi_k\widehat f.
\Ee
Each $P_k$ is bounded on $L^p(\bbR^{d+1})$, $1\le p\le \infty$,
with uniform bounds.

\begin{lemma}\label{pklemmainterpol2d}
Let $d=2$, and suppose that
hypothesis $\cW(w;2) $ holds.
Let
\Be
\label{betastar}
\beta= \beta_*(p,w)=
\frac{3w-13}{6w-20}- \frac{9w-40}{(6w-20)p}
\Ee
and let $r=r(p,w)$ be defined by
\Be \label{rpw}
\frac{1}{r(p,w)}=\frac 12 -\frac{w-2}{6w-20}\Big(3-\frac {10}p\Big)
\Ee
Then,  for $10/3\le p\le w$,
$$\Big\|\sum_k P_k g_k\Big\|_p \le C_\eps \delta^{-\beta-\eps}
\Big\|\Big(\sum_k|g_k|^r\Big)^{1/r}\Big\|_p,
$$
for all families $\{g_k\}$ with $g_k\in \cS(\bbR^{d+1})$.
\end{lemma}

\begin{proof}
By $\cW(w;2)$ and the embedding $L^p(\ell^2)\subset \ell^p(L^p)$
 we have
the inequality
\begin{align}
\Big\|\sum_k P_k g_k \Big\|_{w}
&\le C_\eps \delta^{-(\alpha(w)+\eps)}
\Big(\sum_k\big\|P_k g_k\|_w^w\Big)^{1/w}
\notag
\\
&\lc C_\eps \delta^{-(\alpha(w)+\eps)}
\Big\|\Big(\sum_k|g_k|^w\Big)^{1/w}\Big\|_w.
\label{w-sqfunctionest}
\end{align}

We also observe that for $2\le p\le 4$
\Be \label{l2vect}
\Big\|\Big(\sum_k|P_k g_k|^2\Big)^{1/2}\Big\|_p
\le C (1+\log \delta^{-1} )^{1/2-1/p}
\Big\|\Big(\sum_k|g_k|^2\Big)^{1/2}\Big\|_p
\Ee
Indeed the left hand side is estimated by
\Be \label{weight}
\sup_{\omega\in L^{(p/2)'}} \Big(\sum_k \int |P_k g_k|^2 \omega dx\Big)^{1/2}
\lc
\sup_{\omega\in L^{(p/2)'}} \Big(\sum_k \int |g_k|^2 M_\delta \omega dx\Big)^{1/2}
\Ee
where $M_\delta$ is a Besicovitch-type maximal operator associated to the
light cone which  is bounded on $L^2$ with norm
$O(\sqrt{\log(2+ \delta^{-1})})$ if $\delta<1/2$, see \cite{cord}, \cite{MSS}.
 Thus
H\"older's inequality implies \eqref{l2vect}.

Now we can combine  Proposition \ref{sqlinear} with respect to the double
plates, and $f_k=P_kg_k$,  and
\eqref{l2vect} to obtain
\begin{align}
\Big\|\sum_k P_k g_k \Big\|_{10/3}
&\le C_\eps \delta^{-\frac{1}{20}-\eps}
\Big\|\Big(\sum_k|g_k|^2\Big)^{1/2}\Big\|_{10/3}.
\label{103-sqfunctionest}
\end{align}
After a little arithmetic the claimed bound follows by
interpolation between
\eqref{w-sqfunctionest} and \eqref{103-sqfunctionest}.
\end{proof}

Since $r(p,w)\ge 2$ in Lemma \ref{pklemmainterpol2d}
  we immediately get
\begin{corollary}\label{lemmainterpol2d}
Let $d=2$, suppose that
hypothesis $\cW(w;2) $ holds.
Then  for all  families of functions $\{f_k\}$ with
$\supp \widehat{f_k}\subset \Pi_k^{(\delta)}$
 the estimate \eqref{sqestimate} holds for $10/3\le p\le w$ with
$\beta= \beta_*(p,w)$.
\end{corollary}

In particular note that
$\beta_*(4,w)= \frac{3w-12}{24 w-80}$
so that
$\beta_*(4,6)= 3/32$.  If we use the exponent obtained in Theorem \ref{th1},
{\it i.e.} $w=p_2=190/3$
we get only $\beta_*(4,p_2)=89/720$ which is worse than
$5/44$ exponent that is already known from
 \cite{TV}, \cite{W2}.

For large values of $w$ one can improve on the  result of Corollary
\ref{lemmainterpol2d}.
Our approach will be similar to the one by Tao and Vargas \cite{TV}
in $2+1$ dimensions. By using
$\cW(w;2)$ in that approach one can slightly improve on the previously
known exponents.

\begin{theorem}\label{thmsqimpr}
Let $2\le d\le 4$, and $p$ be as in
\eqref{threecases}.
If hypothesis $\cW(w;d)$ holds then for all  families
 of Schwartz functions
$\{f_k\}$ with
$\supp \widehat{f_k}\subset \Pi_k^{(\delta)}$
 the estimate \eqref{sqestimate} holds with
\begin{equation}
\label{defbetap}
\beta=\,
\mu(p)\, -
\frac{d-1}2 \Big(
\frac{
\frac{d+1}{2(d+3)}-\frac 1p   }
{\frac{d+1}{2(d+3)}+\frac 1p   -\frac{2(p-1)}{(w-1)p}  }
\Big)
 \Big(\frac 1p-\frac{d-2}{2(d-1)}\Big)
\end{equation}




\end{theorem}
The proof (of a slightly more general result)  will be given in \S\ref{five}.


In $2+1$ dimensions Theorem \ref{thmsqimpr} yields inequality \eqref{sqestimate}
for the range $10/3\le p\le w$
 with $\beta$ equal to
\Be\label{betastarstar}
\beta_{**}(p,w)=\frac{1}{2p}  \cdot
\frac{(3p^2-2p-20)w-23p^2+82p-40}{(10+3p)w-23 p+10};
\Ee
in particular we have
$\beta_{**}(4,w)
=\frac{5w-20}{44w-164}$ which (with $p_2\equiv w$) occurs in
Theorem \ref{544thm}.
We compare with \eqref{betastar}. Notice that
$3/32=\beta_*(4,6)<\beta_{**}(4,6)=1/10$.
 A straightforward computation shows
the inequality $\beta_{**}(p,w)<\beta_{*}(p,w)$ holds if and only if
$(9p-30)w^2+(-9p^2-39 p+230)w+23 p(3p-10)>0$ and after factoring we see
that for $10/3<p< w$ we have
$\beta_{**}(p,w)<\beta_{*}(p,w)$ if and only if
$(p-\frac {10}3)(w-\frac{23}{3}) (w-p)>0.$
Thus for any $p\in (10/3,w)$ we have
\Be \beta_{**}(p,w)<\beta_{*}(p,w) \quad \iff w>\frac {23}3
\label{starversusstarstar}
\Ee
 so that the $L^p$ result in
Theorem \ref{thmsqimpr}  is better than the result of Corollary
 \ref{lemmainterpol2d}
in the range
$w>23/3$.
We obtain the following  corollary which yields Theorem \ref{544thm}.

%
%
%
%


\begin{corollary}\label{cone2d}
Let $d=2$ and  suppose that $\cW(w;2)$ holds for some $w>6$. Let
  $10/3<p\le 4$ and let $\alpha> \min \{\beta_*(p,w), \beta_{**}(p,w)\}$
(i.e. $\alpha > \beta_{**}(p,w)\}$ if $w>23/3$).

\noindent Then

(i) the smoothing inequality \eqref{lsm} holds true and

 (ii) the Fourier  multiplier $m_\alpha$ in
 \eqref{conemultiplier} defines a bounded operator on
$L^p(\bbR^3)$.


\end{corollary}
We also observe that by interpolation
we obtain the analogous boundedness results for
the range  $4\le p\le w$  under the assumption that
$\alpha> \frac 12-\frac 2p+ \frac{4(w-p)}{p(w-4)}
\min\{\beta_*(4,w)),  \beta_{**}(4,w)\}$.

If we use the result of Theorem \ref{th1} in $2+1$ dimensions
({\it i.e.}  hypothesis $\cW(w;2)$  with $w=p_2=190/3$)
we obtain this result for
$\alpha >\beta_{**}(p,\tfrac{190}{3})=
 \frac{501 p^2-134 p-3920}
{2p(501 p+1930)}.
$
which equals $445/3934$ if $p=4$. This represents a slight improvement over
the Tao-Vargas result \cite{TV}  which yields the $L^4$ boundedness
 for $\alpha> \frac{5}{44} =
0.113\overline{63}$; note that
$\frac{445}{3934}\approx
0.11311642...$. We also see from  Corollary \ref{lemmainterpol2d}  that
the validity of \eqref{Wp} for the optimal (conjectured) range $p\ge 6$
implies the $L^4$ boundedness for $\alpha>3/32=0.09375$;
however
it has been conjectured that it should hold for all
$\alpha>0$.

\begin{proof}[Proof of Corollary \ref{cone2d}]
It remains to estimate the $L^p$ norm of the square function. For
part (ii) this is done as in \cite{Mock}, namely one first uses a
weighted $L^2$ bound as in \eqref{weight} together with the
optimal $L^{(p/2)'}$ bound for a Besicovich maximal function
associated with the light cone. Now let $S_k$ be the  region in
$\bbR^2$ obtained by projecting the plate $\Pi^{(\delta)}_k$ to
the $\xi_1-\xi_2$ plane. Now define an operator $\fS_k$ by
$\widehat{\fS_k g}(\tau,\xi) = \eta_k(\xi) \widehat g(\tau,\xi)$
where $\xi=(\xi_1,\xi_2)$ and
 $\eta_k$ is a function adapted to the double of $S_k$, with the property that
 $P_k \fS_k=P_k$. We then obtain
\Be \label{l2vectSk} \Big\|\Big(\sum_k|P_k g|^2\Big)^{1/2}\Big\|_p
\le C (1+\log \delta^{-1} )^{1/2-1/p} \Big\|\Big(\sum_k|\fS_k
g|^2\Big)^{1/2}\Big\|_p \Ee and by C\'ordoba's estimate for a
sectorial square function (\cite{cord}) one  dominates the latter
$L^p$ norm  by $C(\log \delta)^{_C}\|g\|_p$. For part (i) one
argues similarly, except that now one has to use a result for a
Besicovitch maximal function which sends  functions on $\bbR^3$ to
functions on $\bbR^2$;  this variant and its application is
discussed in \cite{MSS}.
\end{proof}

\medskip

\section{ Proof of Theorem \ref{thmsqimpr} }
\label{five}

We work with the operators $P_k$ in \eqref{pkdef} which localize
in Fourier space to the doubles of the plates $\Pi_k^{(\delta)}$.
It will be convenient with the following mixed norm variant of the
\lq Wolff hypothesis'.
\medskip

\noindent {\bf Hypothesis $\boldsymbol{\cW(r,s;d)}$.} \emph{Given
$r\geq2(d+1)/(d-1)$ and $1\leq s\leq r$, we say that hypothesis
$\cW(r,s;d)$ holds if for all $\dt<1$, $\e>0$, and all families of
Schwartz functions $\{h_k\}$, we have
\begin{equation}\label{lrestimate}
\Big\|\sum_k P_k h_k \Big\|_{r}\le C_\eps
 \delta^{-\alpha(r,s)-\eps}
\Big(\sum_k \|h_k\|_r^s\Big)^{1/s},
\end{equation}
where $\alpha(r,s)=\frac{d-1}{2s'}-\frac{d+1}{2r}$.}

We shall prove the following variant of Theorem
\ref{thmsqimpr}.

\begin{theorem}\label{thmsqimprwithrs}
Let $2\le d\le 4$, and $p$ be as in
\eqref{threecases}. Let $q=\frac{2(d+3)}{d+1}$.
If hypothesis $\cW(r,p;d)$ holds then for all  families
 of Schwartz functions
$\{f_k\}$ with
$\supp \widehat{f_k}\subset \Pi_k^{(\delta)}$
 the estimate \eqref{sqestimate} holds with
\begin{equation}
\label{defgammastarp} \beta=\, \mu(p)\, -
\tfrac{d-1}2 \Bigl[ \frac{ \frac1q
-\frac 1p   } {\frac1q
+\frac 1p -\frac2r} \Bigr]
 \Big(\tfrac 1p-\tfrac{d-2}{2(d-1)}\Big).
\end{equation}
\end{theorem}

Theorem \ref{thmsqimpr} is an immediate consequence of
Theorem \ref{thmsqimprwithrs}, by the following observation.


\begin{lemma}\label{lrwolff}
Let $w\geq\frac{2(d+1)}{d-1}$ and fix $p\in [2,w]$.
Then $\cW(w;d)$ implies $\cW(r,p;d)$ with $r=p'(w-1)$.
\end{lemma}

\begin{proof}
This follows by interpolation between the Wolff inequality
(i.e. \eqref{w-sqfunctionest} in $d$ dimensions) and the trivial bound
$\|\sum_k P_kh_k\|_\infty\lc \sum_k\|h_k\|_\infty$.
\end{proof}

\noindent {\it Remark:}  In \cite{GSS} we establish certain cases of the
mixed norm inequality $\cW(r,s;d)$ which do not simply  follow by
interpolation from the
original Wolff inequality (as formulated in $\cW(w;d)$). In
such cases Theorem  \ref{thmsqimprwithrs} leads to further
improvements of Theorem \ref{544thm}.

\medskip

To establish Theorem
\ref{thmsqimprwithrs} we  shall work with the following

\noindent{\bf Hypothesis $\boldsymbol{{\mathcal{S\!Q}}(\ga,p)}$.}
 \emph{For all $\delta<1$, $\eps>0$}
\begin{equation}\label{hypothesisga}
\Big\|\sum_k h_k \Big\|_p\le C_\eps \delta^{-\gamma-\eps}
\Big\|\Big(\sum_k|h_k|^2\Big)^{1/2}\Big\|_p,
\end{equation}
\emph{provided that $\supp \widehat{h_k}\subset \Pi_k^{(\delta)}$.}

By Proposition
\ref{sqlinear}
we know  already that for $p>2(d+3)/(d+1)$
 this inequality holds true with the exponent
$\gamma=\mu(p)= d/4-(d+1)/2p$ and we seek an improvement in the ranges
\eqref{threecases}.

We use the hypothesis $\cW(r,p;d)$
 to prove the following proposition, which
amounts to an improved version of Proposition 5.4 in
 \cite{TV} (where
the case $r=\infty$  was considered  in the $2+1$ dimensional situation). As in \S
\ref{two} we work with a covering $\cQ(\delta^{-1/2})$ of
$\sqrt{1/\delta}$ cubes.

\begin{proposition}\label{Lrproposition} Let $d\geq2$, $2<p<r$, and suppose that
hypotheses $\cW(r,p;d)$ and $\mathcal{S\!Q}(\ga,p)$ hold. Then,
for all functions $h_k$ with $\supp \widehat{h_k}\in
\Pi_k^{(\delta)}$ we have
\begin{equation}\label{lpLr}
\Big(\sum_{Q\in \cQ(\delta^{-1/2})} \Big\|\sum_k h_k
\Big\|_{L^r(Q)}^p \Big)^{1/p}\le C_\eps
\delta^{-\frac{\gamma+\alpha(p)}2- \eps}
\Big\|\Big(\sum_k|h_k|^2\Big)^{1/2}\Big\|_{L^p(\SRd)}\ .
\end{equation}
\end{proposition}


%
%

\begin{proof}
We group the indices $k$ (and therefore the
corresponding plates $\Pi_k^{(\delta)}$)  into $O(\delta^{-(d-1)/4})$ disjoint families $S_l$ so that
$\text{dist}(\omega_k,\omega_{k'})
\lc \delta^{1/4}$ for $k,k'\in S_l$. Define
$$G_l=\sum_{k\in S_l} g_k.$$

As in the proof of Proposition \ref{propbilinear} we also work
with the functions $\psi_Q$ adapted to the cubes $Q\in \cQ(\delta^{-1/2})$.
By the support property of $\widehat {\psi_Q}$  the Fourier transform of
$\psi_QG_l$ is supported in a $C\sqrt{\delta}$ plate and these plates form
an essentially disjoint plate family.
Therefore
\begin{align}
\Big\|\sum_lG_l\Big\|_{L^r(Q)} &\lc
\Big\|\psi_Q \sum_lG_l\Big\|_r
\notag
\\
&\leps \delta^{-\frac{\alpha(r,p)}2 }
\Big(\sum_l\big\|\psi_QG_l\big\|_r^p\Big)^{1/p},
\end{align}
by 
hypothesis $\cW(r,p;d)$ with $\delta$ replaced by $\sqrt{\delta}$.
By the support property of $\widehat{\psi_Q G_l}$ and Young's
inequality \Be \|\psi_QG_l\|_r
 \lc \delta^{\frac{d+1}4(\frac 1p-\frac 1r)}
\|\psi_QG_l\|_p
\label{young}
\Ee
and therefore
$$
\Big(\sum_Q\big\|\sum_lG_l\big\|^p_{L^r(Q)} \Big)^{1/p} \leps
\delta^{-\frac{\alpha(r,p)}2+ \frac{d+1}4(\frac 1p-\frac 1r)}
\Big(\sum_{Q,l}\big\|\psi_QG_l\big\|^p_p \Big)^{1/p}.
$$
A little algebra shows
$$-\frac{\alpha(r,p)}2  +
\frac{d+1}4\Big(\frac 1p-\frac 1r\Big)=-\frac{\alpha(p) }2.$$  From some
straightforward  estimation using the decay of the $\psi_Q$ we
also obtain \Be \Big(\sum_Q\big\|\sum_lG_l\big\|^p_{L^r(Q)}
\Big)^{1/p}\leps \delta^{-\alpha(p)/2}
\Big(\sum_l\|G_l\|_p^p\Big)^{1/p} \Ee

As $\widehat G_l$ is supported in a $C\sqrt\delta$ plate we may
use rescaling arguments as in the proof of Lemma \ref{bil1} to
deduce from the hypothesis $\cS\!\cQ(\gamma,p)$ applied with
parameter $\sqrt \delta$ that
$$\|G_l\|_p \leps \delta^{-\gamma/2} \Big\|\Big(
\sum_{k\in S_l}|g_k|^2\Big)^{1/2}\Big\|_p
$$
and hence
\begin{align*}
\Big(\sum_Q\big\|\sum_lG_l\big\|^p_{L^r(Q)} \Big)^{1/p}
&\le C_\eps
\delta^{-\frac{\alpha(p)+\gamma}{2}-\eps}
\Big(\sum_l\Big\|\Big(
\sum_{k\in S_l}|g_k|^2\Big)^{1/2}\Big\|_p^p\Big)^{1/p}
\\
&\lc C_\eps
\delta^{-\frac{\alpha(p)+\gamma}{2}-\eps}
\Big\|\Big(
\sum_{k}|g_k|^2\Big)^{1/2}\Big\|_p
\end{align*}
which is the assertion.
\end{proof}


\begin{proof}[Proof of Theorem \ref{thmsqimprwithrs}, cont.]
We begin by observing that  Hypothesis
$\cS\!\cQ(\mu(p),p)$ holds by Proposition \ref{sqlinear}.

Assuming that $\cS\!\cQ(\gamma,p)$ holds for some $\gamma\le
\mu(p)$ the following estimate for bilinear expressions is an
immediate consequence of Proposition \ref{Lrproposition}
\Be \Big(\sum_{Q\in \cQ}
    \Bigl\|\,\bigl(\sum_{\om_k\in\O}
f_k\bigr)\,\bigl(\sum_{\om_{k'}\in\O'}
g_{k'}\bigr)\,\Bigr\|_{L^{r/2}(Q)}^{p/2}\Big)^{2/p}\,\leps\,\,
\delta^{-\alpha(p)-\gamma}
\Bigl\|\,\bigl(\sum_{\om_k\in\O}|f_k|^2\bigr)^\frac12\,\Bigr\|_p
\,
\Bigl\|\,\bigl(\sum_{\om_{k'}\in\O'}|g_{k'}|^2\bigr)^\frac12\,\Bigr\|_p.
\label{lpLr-bil}
\Ee

We now assume that $\O$ and $\O'$ are separated as in Proposition
\ref{propbilinear}
and interpolate the inequalities \eqref{lpLr-bil} and \eqref{ss2}
with $q=2(d+3)/(d+1)$.
As a result we obtain

\begin{multline*}
\Big(\sum_{Q\in \cQ}
    \Bigl\|\,\bigl(\sum_{\om_k\in\O}
f_k\bigr)\,\bigl(\sum_{\om_{k'}\in\O'}
g_{k'}\bigr)\,\Bigr\|_{L^{p/2}(Q)}^{p/2}\Big)^{2/p}\,
\leps\,\,
\delta^{-2\Gamma(p,\gamma)}
\Bigl\|\,\bigl(\sum_{\om_k\in\O}|f_k|^2\bigr)^\frac12\,\Bigr\|_p
\,
\Bigl\|\,\bigl(\sum_{\om_{k'}\in\O'}|g_{k'}|^2\bigr)^\frac12\,\Bigr\|_p.
\end{multline*}
%
where
$$
\Gamma(p,\gamma)=
(1-\vartheta)\mu(p)+ \vartheta \frac{\alpha(p)+\gamma}{2}
\quad \text { with }  \vartheta=
\Big(\frac 1q-\frac 1p\Big)\Big/\Big(\frac 1q-\frac 1r\Big).
$$
By Lemma \ref{bil1} we also obtain
\Be
\label{lpLplinear}
    \Bigl\|\sum_{k}
f_k\Bigr\|_{p}\, \leps\,\,
\delta^{-\Gamma(p,\gamma)}
\Bigl\|\bigl(\sum_{\om_k\in\O}|f_k|^2\bigr)^{1/2}\Bigr\|_p.
\Ee

The assumption $p<\frac{2(d-1)}{d-2}$ in \eqref{threecases}
implies that $\alpha(p)<\Gamma(p,\gamma) \le \mu(p)$ provided that
$\alpha(p)<\gamma\le\mu(p)$. Moreover $\gamma= \Gamma(p,\gamma)$
if and only if $\gamma$ equals
$$
\gamma_*
=\frac{1}{1-\vth/2}\big((1-\vth)\mu(p)+\vth \frac{\alpha(p)}{2}\big)\,
=\,\mu(p)- \frac{\vartheta}{2-\vartheta}\big(\mu(p)-\alpha(p)\big).
$$
The fixed point  is contained in the interval $(\alpha(p),\mu(p))$
and one observes  that $\Gamma(p,\gamma)<\gamma$ for
$\gamma^*<\gamma\leq\mu(p)$. Thus, if we define a sequence
$\gamma_n$ by setting $\gamma_0=\mu(p)$ and
$\gamma_{n+1}=\Gamma(p,\gamma_n)$ for $n\ge 0$, then $\gamma_n$ is
decreasing and bounded below and converges to $\gamma^*$. We
compute that $\vth/(2-\vth)=(1/q-1/p)\big /(1/q+1/p-2/r)$ and
$\alpha(p)-\mu(p)= (d-2)/4-(d-1)/2p$ and see that $\gamma_*$ is
equal to the right hand side of \eqref{defgammastarp}.
Thus \eqref{lpLplinear} and an iteration yield
the assertion of the theorem.
\end{proof}

\bibliographystyle{plain}

\small

\end{document}